\documentclass[12pt]{article}
\usepackage{mathrsfs}
\usepackage{stmaryrd}
\usepackage{amsfonts,amsmath,amssymb,amscd}
\usepackage{shadow}
\usepackage{graphicx}
\usepackage{color}
\usepackage{longtable}
\allowdisplaybreaks

\newtheorem{theo}{Theorem}[section]

\newtheorem{lemm}[theo]{Lemma}
\newtheorem{prop}[theo]{Proposition}
\newtheorem{coro}[theo]{Corollary}

\def\qed{\hfill \rule{4pt}{7pt}}

\def\pf{\noindent {\it{Proof.} \hskip 2pt}}

\numberwithin{equation}{section}

\makeatletter

\newcommand{\Rmnum}[1]{\expandafter\@slowromancap\romannumeral #1@}
\makeatother

\begin{document}
\begin{center}
{\large\bf The Sorting Index  and Permutation Codes}
\end{center}

\begin{center}
William Y.C. Chen$^1$, George Z. Gong$^2$, Jeremy J.F. Guo$^3$

Center for Combinatorics, LPMC-TJKLC\\
Nankai University, Tianjin 300071, P.R. China

$^1$chen@nankai.edu.cn, $^2$sxgong@mail.nankai.edu.cn, $^3$guojf@mail.nankai.edu.cn 

\end{center}

\begin{abstract}
In the combinatorial study of the coefficients of a bivariate polynomial that generalizes both the length and the reflection length generating functions for finite Coxeter groups, Petersen introduced a new Mahonian statistic $sor$, called the sorting index.
Petersen proved  that the pairs of statistics
$(sor,cyc)$ and   $(inv,rl\textrm{-}min)$ have the
same joint distribution over the symmetric group,
and asked
 for a combinatorial proof  of this fact. In answer to the
 question of Petersen, we observe a
 connection between the sorting index and the B-code of a permutation defined by
 Foata and Han, and we show that the bijection of Foata and Han serves the
 purpose of mapping    $(inv,rl\textrm{-}min)$ to $(sor,cyc)$.
We also give a type $B$ analogue of the Foata-Han bijection,
 and we derive the equidistribution of
 $(inv_B,{\rm Lmap_B},{\rm Rmil_B})$ and
$(sor_B,{\rm Lmap_B},{\rm Cyc_B})$ over signed permutations.
So we get a combinatorial interpretation of
Petersen's equidistribution of $(inv_B,nmin_B)$ and $(sor_B,l_B')$.
Moreover, we show that the six pairs of set-valued statistics
 $\rm (Cyc_B,Rmil_B)$, $\rm(Cyc_B,Lmap_B)$, $\rm(Rmil_B,Lmap_B)$, $\rm(Lmap_B,Rmil_B)$, $\rm(Lmap_B,Cyc_B)$ and
$\rm(Rmil_B,Cyc_B)$  are equidistributed over signed permutations.
For Coxeter groups of type $D$,  Petersen showed that
 the two statistics $inv_D$ and $sor_D$
are equidistributed.  We introduce two statistics
 $nmin_D$ and $\tilde{l}_D'$ for elements of $D_n$ and we prove that the two pairs  of statistics $(inv_D,nmin_D)$ and $(sor_D,\tilde{l}_D')$
 are equidistributed.
\end{abstract}



\section{Introduction}

This paper is concerned with a combinatorial study of the Mahonian statistic
$sor$, introduced by Petersen \cite{Petersen}. This statistic is also interpreted by Wilson \cite{Wilson1,Wilson2}  as the total distance moved rightward in the random generation of a permutation based on the Fisher-Yates shuffle algorithm.

Let $[n]=\{1,2 ,\ldots, n\}$. The set of permutations of $[n]$ is denoted by $S_n$.
 Let us recall the definition of the sorting index of a permutation $\sigma$ in $S_n$.
Notice that $\sigma$  has a unique decomposition into transpositions
 $$\sigma=(i_1,j_1)(i_2,j_2)\cdots(i_k,j_k)$$ such that
$$j_1<j_2<\cdots<j_k$$ and $$i_1<j_1,i_2<j_2,\ldots,i_k<j_k.$$
The sorting index is defined by
$$sor(\sigma)=\sum_{r=1}^k(j_r-i_r).$$

 Based on the cycle decomposition of a permutation, Foata and Han \cite{FH2} introduced the  B-code of a permutation. We observe that
  the sorting index of a permutation can be easily expressed in terms of
  its B-code. Given a permutation $\sigma\in S_n$ with B-code $b=(b_1,b_2,\ldots,b_n)$, it can be seen that the sorting index of $\sigma$
  is given by
   $$sor(\sigma)=\sum^n_{i=1}(i-b_i).$$

Petersen \cite{Petersen} has shown that the sorting index $sor$ is a
Mahonian statistic, that is, it has the same distribution as the
number of inversions. He also introduced the sorting indices for
Coxeter groups of type $B$ and type $D$ and showed that they are Mahonian as well.

Let us recall some notation and terminology. For $n
\geq 1$, given a permutation $\sigma=\sigma_1\sigma_2\cdots\sigma_n\in S_n$, a pair
$(\sigma_i,\sigma_j)$ is called an inversion if $i<j$ and
$\sigma_i>\sigma_j$. Let $inv(\sigma)$ denote the number of
inversions of $\sigma$.  An element  $\sigma_i$  is said to be a
right-to-left minimum of $\sigma$  if $\sigma_i<\sigma_j$ for all $j>i$.
The number of right-to-left minima of $\sigma$ is denoted by
$rl\textrm{-}min(\sigma)$. The number of elements of $\sigma$ that are not right-to-left minima is denoted by $nmin(\sigma)$.
Similarly, one can define
a left-to-right maximum. The number of left-to-right maxima of $\sigma$ is denoted by $lr\textrm{-}max(\sigma)$.  The number of cycles of $\sigma$ is denoted by
$cyc(\sigma)$. The reflection
length of $\sigma$, denoted $l'(\sigma)$,  is the minimal number of
transpositions needed to express $\sigma$.

By using two factorizations of the diagonal sum, i.e.,~$\sum_{\sigma\in S_n}\sigma $, in the group algebra
$\mathbb{Z}[S_n]$,
Petersen has shown that $(sor,cyc)$ and $(inv,rl\textrm{-}min)$ have
the same joint distribution by deriving the following generating function formulas:
$$\sum_{\sigma\in
S_n}q^{sor(\sigma)}t^{cyc(\sigma)}=\sum_{\sigma\in
S_n}q^{inv(\sigma)}t^{rl\textrm{-}min(\sigma)}=t(t+q)\cdots(t+q+q^2+\cdots+q^{n-1}). $$
He raised the question of finding
 a  bijection that maps a permutation with inversion number
$k$ to a  permutation with sorting index $k$.
We find that a bijection constructed by Foata and Han \cite{FH2}
 on $S_n$ serves the purpose of mapping
 $(inv,rl\textrm{-}min)$ to
$(sor,cyc)$.

The bijection of Foata and Han is devised for the purpose of
deriving the equidistribution of  the six pairs of set-valued
statistics $\rm (Cyc,Rmil)$, $\rm(Cyc,Lmap)$, $\rm(Rmil,Lmap)$,
$\rm(Lmap,Rmil)$, $\rm(Lmap,Cyc)$ and $\rm(Rmil,Cyc)$ over $S_n$.
It should be mentioned that the equidistribution of the
 three pairs of set-valued statistics $\rm(Lmap,Cyc)$, $\rm
 (Cyc,Lmap)$, $\rm(Lmap,Rmil)$
reduces to  the equidistribution of the three pairs
of integer-valued statistics $(lr\text{-}max,cyc)$, $(cyc,lr\text{-}max)$
  and $(lr\text{-}max,lr\text{-}min)$
  established by  Cori
\cite{Cori} by employing labeled Dyck paths and the Ossona de Mendez
Rosenstiehl  algorithm \cite{OMR}
 on hypermaps.

As for Coxeter groups of type $B$, the sorting index
can be analogously defined and it is Mahonian, see  Petersen \cite{Petersen}.
Let $sor_B, inv_B, nmin_B$ and $l_B'$ denote
the statistics on signed permutations analogous to $sor, inv, nmin$ and $l'$ for permutations.
Petersen obtained
the following formulas for the joint distributions of $(inv_B,nmin_B)$
and $(sor_B,l_B')$:
$$\sum_{\sigma\in
B_n}q^{sor_B(\sigma)}t^{l_B^{'}(\sigma)}=\sum_{\sigma\in
B_n}q^{inv_B(\sigma)}t^{nmin_B(\sigma)}= \prod_{i=1}^n (1+t[2i]_q-t) . $$

 We shall present a bijection on $B_n$ which implies the equidistribution of
 $(inv_B,{\rm Lmap_B},$\\${\rm Rmil_B})$ and
$(sor_B,{\rm Lmap_B},{\rm Cyc_B})$ where $\rm Lmap_B$, $\rm Rmil_B$ and $\rm Cyc_B$ are  set-valued statistics.
 In particular, this bijection transforms $(inv_B,nmin_B)$  to  $(sor_B,l_B')$.
 We introduce the A-code and the B-code of a signed permutation,
  which are analogous to the A-code and the B-code of a permutation.
   We  show that  the triple of
     statistics $(inv_B,{\rm Lmap_B},{\rm Rmil_B})$ of a signed permutation can be computed from its A-code  whereas the triple of statistics $(sor_B,{\rm Lmap_B},{\rm Cyc_B})$  can be computed from its B-code.
To be more specific, let  $\sigma$ be a signed permutation in $B_n$ with   A-code  $c$.
 Let $\sigma'$ be the signed permutation in $B_n$ with   B-code $c$.
  Then the triple of
     statistics $(inv_B,{\rm Lmap_B},{\rm Rmil_B})$ of $\sigma$
  coincides with the triple of statistics $(sor_B,{\rm Lmap_B},{\rm Cyc_B})$  of $\sigma'$.
We also show that the six pairs of set-valued statistics $\rm
(Cyc_B,Rmil_B)$, $\rm(Cyc_B,Lmap_B)$, $\rm(Rmil_B,Lmap_B)$,
$\rm(Lmap_B,Rmil_B)$, $\rm(Lmap_B,Cyc_B)$ and $\rm(Rmil_B,Cyc_B)$
are equidistributed over $B_n$. As a consequence, we see that
 the four pairs of statistics $(sor_B,l'_B)$, $(inv_B,nmin_B)$,
$(inv_B,nmax_B)$ and $(sor_B,nmax_B)$ are equidistributed over $B_n$.

For Coxeter groups of type $D$, let $sor_D$ and $inv_D$ denote the statistics analogous  to $sor$ and $inv$. Let $D_n$ denote the subgroup of $B_n$ consisting of all signed permutations with an even number of minus signs.
In this case, Petersen has shown that $sor_D$ and  $inv_D$ have the same
generating function, that is,
$$
\sum_{\sigma\in
  D_n}q^{sor_D(\sigma)}=\sum_{\sigma\in
  D_n}q^{inv_D(\sigma)}
  ={[n]}_q\prod_{r=1}^{n-1}{[2r]}_q. $$

We shall introduce two statistics $nmin_D$ and $\tilde{l}_D'$ analogous to
 $nmin$ and $l'$, and we shall construct a bijection in order to show that
  the pairs of statistics $(inv_D,nmin_D)$ and $(sor_D,\tilde{l}_D')$ are equidistributed over $D_n$.
  Moreover, we prove that the bivariate  generating functions
 for $(inv_D,nmin_D)$ and $(sor_D,\tilde{l}_D')$
 are both equal to
 $$D_n(q,t)=\prod_{r=1}^{n-1}(1+q^rt+qt\cdot{[2r]}_q)~.$$

\section{The bijection of Foata and Han}

In this section, we give a brief description of Foata and Han's
bijection \cite{FH2} on permutations.
 Then we shall show that this bijection
indeed transforms $(inv,rl\textrm{-}min)$ to $(sor,cyc)$.

The  group of permutations of $[n]$ is also known as a Coxeter group
of type $A$. The \emph{length} of a permutation $\sigma\in S_n$,
denoted by $l(\sigma)$, is defined to be the minimal number of
adjacent transpositions needed to express $\sigma$. It is not
difficult to see that $inv(\sigma)=l(\sigma)$.

We  adopt the notation of Foata and Han \cite{FH2}. They have
investigated several set-valued statistics which are defined as
follows. Given a permutation $\sigma\in S_n$, it can be decomposed
as a product of disjoint cycles whose minimum elements are
$c_1,c_2,\ldots,c_r$. Define $\textrm{Cyc} ~\sigma$ to be the set
$$\textrm{Cyc}~\sigma=\{c_1,c_2,\ldots,c_r \}.$$
Let $\omega=x_1x_2\cdots x_n$ be a word in which the letters are
positive integers. The {left to right maximum place set} of $\omega$,
denoted by $\textrm{Lmap}~\omega$, is the set of all places $i$ such
that $x_j<x_i$ for all $j<i$, while the {right to left minimum letter
set} of $\omega$, denoted by $\textrm{Rmil}~\omega$,  is the set of
all letters $x_i$ such that $x_j>x_i$ for all $j>i$. For a permutation
 $\sigma$ of $[n]$, recall that
$lr\textrm{-}max(\sigma)$ is the number of left-to-right maxima of $\sigma$,
$rl\text{-}min(\sigma)$ is the number of right-to-left minima of $\sigma$, and $cyc(\sigma)$ is the number of cycles of $\sigma$.
 It is easy to see that the cardinalities of
$\textrm{Lmap}~\sigma$, $\textrm{Rmil}~\sigma$ and
$\textrm{Cyc}~\sigma$ reduce to $lr\text{-}max (\sigma)$,
$rl\text{-}min(\sigma)$ and $cyc(\sigma)$, respectively.

 The Lehmer
code \cite{Lehmer} of a permutation
$\sigma=\sigma_1\sigma_2\cdots\sigma_n$ of $[n]$ is defined to be the
sequence $\textrm{Leh}~\sigma=(a_1,a_2,\ldots,a_n)$, where
$$a_i=|\{j\mid 1\leq j\leq i,\sigma_j\leq \sigma_i\}|.$$
 Let $\textrm{SE}_n$ denote the set of integer sequences
 $(a_1,a_2,\ldots,a_n)$ such that $1\leq a_i \leq i$ for all $i$.
Then ${\rm Leh}\colon S_n\longrightarrow \textrm{SE}_n$ is a bijection. Foata and Han
\cite{FH2} defined the A-code of a permutation $\sigma$ to be
a sequence $$\textrm{A-code}~\sigma=\textrm{Leh}~ \textbf{i}\sigma  $$
where $\textbf{i}\colon \sigma\mapsto\sigma^{-1}$ denotes the inverse operation on $S_n$ with respect
to product of permutations. For example, for $\sigma=3\,1\,5\,2\,4$, then $\textbf{i}\sigma=2\,4\,1\,5\,3$. Here a permutation $\sigma=\sigma_1\sigma_2\cdots\sigma_n\in S_n$ standards for a
 one-to-one function on $[n]$ which maps $i$ to $\sigma_i$ for $1\leq i\leq n$. We multiply permutations from right to left, that is, for $\pi,\sigma\in S_n$, we have $\pi\sigma(i)=\pi(\sigma(i))$ for $1\leq i \leq n$.

For an
integer sequence $a=(a_1,a_2,\ldots,a_n)\in\textrm{SE}_n$, define $\textrm{Max}~a$ to be the set $\{i\mid a_i=i  \}$.
Given a permutation $\sigma\in S_n$, Foata and Han \cite{FH2} have shown that
 the A-code leads to a bijection from $S_n$ to $\textrm{SE}_n $ and
 the two set-valued statistics $\textrm{Rmil}$ and $\textrm{Lmap}$ of $\sigma$ are determined by  its A-code, that is,
\begin{equation} \label{R}
  \textrm{Rmil}~\sigma=\textrm{Max}~(\textrm{A-code}~\sigma),
\end{equation}
\begin{equation} \label{L}
  \textrm{Lmap}~\sigma=\textrm{Rmil}~(\textrm{A-code}~\sigma).
\end{equation}
Following the notation in \cite{FH2},
 we rewrite (\ref{R}) and (\ref{L}) as
\begin{equation}\label{a1}
(\textrm{Rmil},\textrm{Lmap})~\sigma=(\textrm{Max},\textrm{Rmil})~\textrm{A-code}~\sigma. \end{equation}

Given a permutation
$\sigma=\sigma_1\sigma_2\cdots\sigma_n\in S_n$, the B-code
can be defined as follows. For $1\leq i\leq n$, let $k_i$ be
the smallest integer $k\geq 1$ such that $(\sigma^{-k})(i)\leq i$,
where $\sigma$ is considered as a bijective function on $[n]$.
Then $b_i=(\sigma^{-k_i})(i)$. In fact, the B-code of a permutation
can be easily
determined by the cycle decomposition.
To compute $b_i$, we assume that $i$ appears in a cycle $C$.
If $i$ is the smallest element of $C$, then we set $b_i=i$. Otherwise,
we choose $b_i$ to be the element $j$ of $C$
such that $j<i$ and $j$ is the closest to $i$. Notice that $C$
 is viewed as a directed cycle and the distance from $j$ to $i$ is
 meant to be the number of steps to reach $i$ from $j$ along the cycle.
For example, let $\sigma=2\,4\,5\,1\,3$. Using the
cycle decomposition $\sigma=(1~2~4)(3~5)$, we get the B-code
$(1,1,3,2,3)$.

Foata and Han have shown that the B-code is a bijection from $S_n$ to $\textrm{SE}_n $ and the pair of set-valued statistics $\rm(Cyc,Lmap)$ of $\sigma$ can be determined by the  B-code of $\sigma$, that is,
\begin{equation}\label{a2}
\expandafter{\rm(Cyc,Lmap)}~\sigma=\expandafter{\rm(Max,Rmil)~B\text{-}code}~ \sigma.
\end{equation}

Combining the A-code and the B-code, Foata and Han \cite{FH2}
found a bijection $\phi$ on $S_n$ as given by
 $$\phi=\rm (B\text{-}code)^{-1}\circ A\text{-}code.$$
The bijection $\phi$  implies the following
equidistributions.

\begin{theo}[Foata and Han \cite{FH2}]\label{thm-set-A}
The six pairs of set-valued statistics  $\rm (Cyc,Rmil)$,
$\rm(Cyc,Lmap)$, $\rm(Rmil,Lmap)$, $\rm(Lmap,Rmil)$,
$\rm(Lmap,Cyc)$, $\rm(Rmil,Cyc)$ are equidistributed over $S_n$:
\end{theo}
\begin{center}
 \begin{tabular}{ccccccccccc}
 $S_n$ &  $\stackrel{\textbf{i}}{\longrightarrow}$ &  $S_n$   &  $\stackrel{\phi^{-1}}{\longrightarrow}$    & $ S_n $ & $ \stackrel{\textbf{i}}{\longrightarrow}$ & $ S_n $ & $\stackrel{\phi }{\longrightarrow}$ & $ S_n $ &
$\stackrel{\textbf{i} }{\longrightarrow}$ &$ S_n $ \\[5pt]
    $\rm Cyc\choose Rmil $   &         & $\rm Cyc \choose Lmap $   &   & $\rm Rmil\choose Lmap $  &      & $\rm Lmap\choose Rmil  $ & & $\rm Lmap\choose Cyc$ &    &$\rm Rmil\choose Cyc$.
 \end{tabular}
\end{center}

We now turn to the sorting index.
Petersen has shown that the pairs of statistics $(sor,cyc)$ and
$(inv,rl\textrm{-}min)$ have the same joint distribution over
permutations and asked for a combinatorial interpretation of
 this fact. We shall show that the map $\phi$
  transforms the pair of statistics
  $(inv,rl\textrm{-}min)$ of a permutation $\sigma$
  to the pair of statistics $(sor,cyc)$ of the permutation $\phi(\sigma)$.
The following lemma shows that the pair of statistics
$(inv,rl\textrm{-}min)$ of $\sigma$ can be computed from the A-code of $\sigma$.

\begin{lemm}\label{a-typeA}
Let $\sigma$ be a permutation in $S_n$ with {\rm A-code} $a=(a_1,a_2,\ldots,a_n)$.
Then we have
\begin{equation}\label{inv-a}
  inv(\sigma)=\sum^n_{i=1}(i-a_i),
\end{equation}
  and
  \begin{equation}\label{rl-a}
   rl\textrm{-}min(\sigma)=|\textrm{Max}~a|.
   \end{equation}
       \end{lemm}

\pf By the definition of the A-code, we find
$$inv(\sigma)={n\choose 2}-\sum^n_{i=1}(a_i-1),$$
which can be rewritten as
$$\sum^n_{i=1}(i-a_i).$$ From (\ref{a1}) it follows that
$rl\textrm{-}min(\sigma)=|\textrm{Rmil}~\sigma|=|\textrm{Max}~a|$,
as desired. \qed

The following lemma shows that the pair of statistics
$(sor,cyc)$ of $\sigma$ can be recovered from the B-code.

\begin{lemm}\label{b-typeA}
Let $\sigma$ be a permutation in $S_n$ with {\rm B-code} $b=(b_1,b_2,\ldots,b_n)$.
Then we have
  \begin{equation}\label{sor-a}
  sor(\sigma)=\sum^n_{i=1}(i-b_i),
  \end{equation}
  and
  \begin{equation}\label{cyc-a}
  cyc(\sigma)=|\textrm{Max}~b|.
  \end{equation}
\end{lemm}

\pf  Let us examine the algorithm of Foata and Han  to recover a permutation $\sigma$ from its B-code
   $b=(b_1,b_2,\ldots,b_n)\in \textrm{SE}_n$.
  Start with the identity permutation $\sigma^{(0)}=12\cdots n$. For $1\leq i\leq n$,
the permutation $\sigma^{(i)}$ is obtained by exchanging $i$ and the letter at the $b_i$-th place in $\sigma^{(i-1)}$. Notice that it may happen that $i=b_i$. Then the resulting permutation $\sigma^{(n)}$ is precisely the permutation with B-code $b$, that is, $\sigma=\sigma^{(n)}$.
So we may write
  $\sigma^{(i)}=\sigma^{(i-1)}(b_i,i)$, where $(b_i,i)$ is  called
   a transposition even when $b_i=i$. Thus we obtain a decomposition of $\sigma$ into transpositions
$$\sigma=(b_1,1)(b_2,2)\cdots(b_n,n).$$
Then by the definition of the sorting index, we see that
 $$sor(\sigma)=\sum^n_{i=1}(i-b_i).$$
It follows from (\ref{a2}) that $cyc(\sigma)=|\textrm{Cyc}~\sigma|=|\textrm{Max}~b|$. This completes the proof.\qed

Combining Lemma \ref{a-typeA} and Lemma \ref{b-typeA}, we conclude that the bijection $\phi=\rm (B\text{-}code)^{-1}\circ A\text{-}code$
transforms  $(inv,rl\textrm{-}min)$ to $(sor,cyc)$, that is, for any $\sigma\in S_n$,
 $$(inv,rl\textrm{-}min)~\sigma=(sor,cyc)~\phi(\sigma).$$

By Theorem \ref{thm-set-A}, the bijection $\phi$ preserves the set-valued statistic $\textrm{Lmap}$. Since
\[ lr\textrm{-}max (\sigma)=| \textrm{Lmap}~\sigma|,\]
$\phi$ preserves the statistic $lr\textrm{-}max$.
Observing that
\[ rl\textrm{-}min(\sigma)=lr\textrm{-}max(\textbf{i}\sigma),\]
 we arrive at the following equidistributions.

\begin{theo}
The four pairs of statistics $(sor,cyc)$, $(inv,rl\textrm{-}min)$,
$(inv,lr\textrm{-}max)$ and $(sor,lr\textrm{-}max)$ are
equidistributed over $S_n$:
 \end{theo}
\begin{center}
 \begin{tabular}{ccccccc}
 $S_n$ &  $\stackrel{\phi^{-1}}{\longrightarrow}$ &  $S_n$   &  $\stackrel{\textbf{i} }{\longrightarrow}$    & $ S_n $ & $ \stackrel{\phi}{\longrightarrow}$ & $ S_n $\\[5pt]
   $sor\choose cyc $   &         & $ inv\choose rl\textrm{-}min $   &   & $inv\choose lr\textrm{-}max $  &      & $sor\choose lr\textrm{-}max  $.                      \\[5pt]
 \end{tabular}
\end{center}

\section{A bijection on signed permutations}

In this section, we construct a bijection which serves as
a combinatorial interpretation of  the
equidistribution of the pairs of statistics $(inv_B,nmin_B)$ and
$(sor_B,l_B^{'})$ over signed permutations. In fact, this bijection implies the equidistribution of $(inv_B,{\rm Lmap_B},{\rm Rmil_B})$ and
$(sor_B,{\rm Lmap_B},{\rm Cyc_B})$ over $B_n$. Moreover,
we show that the six pairs
of set-valued statistics $\rm (Cyc_B,Rmil_B)$, $\rm(Cyc_B,Lmap_B)$,
 $\rm(Rmil_B,Lmap_B)$, $\rm(Lmap_B,Rmil_B)$,
$\rm(Lmap_B,Cyc_B)$ and $\rm(Rmil_B,Cyc_B)$  are equidistributed over
$B_n$.

Let us recall some definitions.
 The hyperoctahedral group $B_n$ is the group of  bijections $\sigma$ on $\{1,2,\ldots,n,\bar{1},\bar{2},\ldots,\bar{n}\}$
 such that $\sigma(\bar{i})=\overline{\sigma(i)}$ for $i=1,2,\ldots,n$, where $\bar{i}$ denotes $-i$.
Clearly, one can represent an element $\sigma\in B_n$ by a signed permutation $a_1a_2\cdots a_n$ of $[n]$, that is,
a permutation of $[n]$ with some elements associated with the minus sign.

The group  $B_n$ has the following Coxeter generators
$$S^B=\{(\bar{1},1),(1,2),(2,3),\ldots,(n-1,n)\}.$$
The set of   reflections of $B_n$ is
$$T^B=\{(i,j):1\leq i<j\leq n\}\cup\{(\bar{i},j):1\leq i\leq j\leq n\},  $$
where the transposition $(i,j)$ means to exchange   $i$ and $j$ and exchange
$\bar{i}$ with $\bar{j}$ provided that $i\neq \bar{j}$, and $(\bar{i},i)$ means to exchange $i$ and $\bar{i}$. For $\sigma\in B_n$, let
 $N(\sigma)$ denote the number of negative elements in the
 signed permutation notation.

As for permutations, Petersen \cite{Petersen}   defined
the sorting index for a singed permutation.
Let  $\sigma$ be a signed permutation in $B_n$. By using the straight selection sort algorithm \cite{knuth} of type $B$, Petersen has shown that
$\sigma$ has a unique
factorization into  a product of signed transpositions in
$T^B$:
\begin{equation}\label{fac-sigma-B}
  \sigma=(i_1,j_1)(i_2,j_2)\cdots(i_m,j_m),
\end{equation}
where $0<j_1<j_2<\cdots<j_m\leq n$.
Then the sorting index of $\sigma$ is defined by
$$sor_B(\sigma)=\sum_{r=1}^m(j_r-i_r-\chi(i_r<0)).  $$
 For example, let
$\sigma=5\,\bar{4}\,\bar{3}\,1\,\bar{2}$.
Then we have $$\sigma=(\bar{1},2)(\bar{3},3)(\bar{2},4)(1,5)$$
and $sor_B(\sigma)=2-(-1)-1+3-(-3)-1+4-(-2)-1+5-1=16$.

For a signed permutation $\sigma\in B_n$, the length of
$\sigma$, denoted $l_B(\sigma)$,  is defined to be the minimal number of transpositions
in $S^B$ needed to express $\sigma$, see \cite{BB}.
The reflection length of
$\sigma$, denoted $l_B^{'}(\sigma)$, is the minimal number of transpositions
in $T^B$ needed to express $\sigma$. The type $B$ inversion number of $\sigma$, denoted  $inv_B(\sigma)$, also denoted $finv$ in \cite{FH},  is defined as
$$inv_B(\sigma)=|\{(i,j):1\leq i<j\leq n,\sigma_i>\sigma_j\}|+|\{(i,j):1\leq i\leq j\leq n,\overline{\sigma_i}>\sigma_j\}|.$$ Like the case of type $A$, we have
$inv_B(\sigma)=l_B(\sigma)$, see \cite[Section 8.1]{BB}.

Recall that for a permutation $\pi\in S_n$, we have $l'(\pi)=n-cyc(\pi)$.
Similarly, the reflection length of a signed permutation can be determined from its cycle decomposition.
A signed permutation $\sigma$ can be expressed as a product of disjoint signed cycles,
 see,  Brenti \cite{Brenti}, Chen and Stanley \cite{Chen2}.
For example, let $\sigma=\bar{6}\,\bar{7}\,4\,\bar{3}\,5\,1\,\bar{2}$.
Then  $\sigma$ can be written as $\sigma=(1~\bar{6})(5)(\bar{7}~\bar{2})(4~\bar{3})$.
A signed cycle is said to be balanced if it contains an even number of minus signs, see \cite{Chen2}. Let $cyc_B(\sigma)$ denote the number of balanced cycles of $\sigma$. It is not difficult to see that  $l_B'(\sigma)=n-cyc_B(\sigma)$.

 We introduce
 some set-valued statistics for signed permutations which are analogous to those for permutations.
 For a signed permutation $\sigma$,
 let $C_1,C_2,\ldots,C_r$ be the balanced signed cycles of $\sigma$.
 Let $c_i$ be the smallest absolute value of elements of $C_i$. Define $\rm{Cyc_B}$ to be the set
$\{c_1,c_2,\ldots,c_r\}$.

Let $\omega=\omega_1\omega_2\cdots\omega_n$ be a word of
length $n$, where $\omega_i$ is an integer. The left to right
maximum place set of $\omega$, denoted $\textrm{Lmap}_\textrm{B} ~\omega$,
and the right to left minimum letter set of $\omega$, denoted
$\textrm{Rmil}_\textrm{B}~\omega$,  are defined as follows,
$$\textrm{Lmap}_\textrm{B}~\omega=\{i\mid\omega_i>|\omega_j|~ \text{for any}~ j<i \},$$
$$\textrm{Rmil}_\textrm{B}~\omega=\{\omega_i\mid0<\omega_i<|\omega_j|~ \text{for any}~j>i\}.$$
When $\sigma$ is a signed permutation, the cardinality of  $\textrm{Lmap}_\textrm{B}~\sigma$ is denoted by $lr\text{-}max_B (\sigma)$
and  the cardinality of $\textrm{Rmil}_\textrm{B}~\sigma$ is denoted by $rl\text{-}min_B(\sigma)$.
Let
$$nmin_B(\sigma)=|\{i:\sigma_i>|\sigma_j| \rm~for~ some~j>i\}|+N(\sigma)$$
and $$ nmax_B(\sigma)=|\{i:0<\sigma_i<|\sigma_j| \rm~for~ some~j<i\}|+N(\sigma).$$
Evidently, $nmin_B(\sigma)=n-rl\text{-}min_B(\sigma)$ and
$nmax_B(\sigma)=n-lr\text{-}max_B (\sigma)$.

The following theorem is due to Petersen \cite{Petersen}.

\begin{theo}\label{th3.1}
The pairs of statistics $(inv_B,nmin_B)$ and $(sor_B,l_B^{'})$ are
equidistributed over $B_n$:
$$  \sum_{\sigma\in
B_n}q^{sor_B(\sigma)}t^{l_B^{'}(\sigma)}=\sum_{\sigma\in
B_n}q^{inv_B(\sigma)}t^{nmin_B(\sigma)}. $$
\end{theo}

Petersen presented two different factorizations of the diagonal sum
$\sum_{\sigma\in B_n}\sigma $
and showed that the two sides of
the above equation  are both equal to
$$ B_n(q,t)=\prod_{i=1}^n (1+t[2i]_q-t).$$

We shall  construct a bijection
 $\psi\colon B_n\longrightarrow B_n$ which transforms $(inv_B,{\rm Lmap_B},{\rm Rmil_B})$ to
$(sor_B,{\rm Lmap_B},{\rm Cyc_B})$.
  This bijection can be described in terms of two codes, the A-code and the B-code for signed permutations.
For a signed permutation
$\sigma=\sigma_1\sigma_2\cdots\sigma_n\in B_n$,
let $\textbf{i}\colon \sigma\mapsto\sigma^{-1}$ denote
the inverse operation on $B_n$ with respect to product of signed permutations.
We define the Lehmer code of the signed permutation $\sigma$ to be the integer sequence
$\textrm{Leh}~\sigma=(a_1,a_2,\ldots,a_n)$,
where for each $i$,
$$a_i=\textrm{sign}~{\sigma_i}\cdot|\{j\mid 1\leq j\leq i,|\sigma_j|\leq |\sigma_i|  \}|. $$
Then the A-code of a signed permutation $\sigma$ is defined to be an integer sequence
$$\textrm{A}\textrm{-code}~\sigma=\textrm{Leh}~\textbf{i}\sigma. $$

Let $\textrm{SE}^\textrm{B}_n$ be the set of integer sequences $(a_1,a_2,\ldots,a_n)$ such that $a_i\in[-i,i]\setminus\{0\}$.
For an integer sequence $a=(a_1,a_2,\ldots,a_n)\in\textrm{SE}^\textrm{B}_n$,
 $\textrm{Max}~a$ stands for the set $\{i\mid a_i=i  \}$.

The following proposition shows that the two set-valued statistics $\expandafter{\rm Rmil_B}$ and $\expandafter{\rm Lmap_B}$ for a signed permutation $\sigma$ can be recovered from the Lehmer code
of $\sigma$. The proof is straightforward, and hence it is omitted.

\begin{prop}\label{leh-b}
${\rm Leh}\colon B_n\longrightarrow {\rm SE}^{\rm B}_n$ is a bijection and for each $\sigma\in B_n$, we have
\begin{equation}
  \expandafter{\rm Rmil_B ~Leh}~\sigma=\expandafter{\rm Rmil_B}~\sigma,
\end{equation}
and
\begin{equation}
\expandafter{\rm Max~Leh}~\sigma=\expandafter{\rm Lmap_B}~\sigma.
\end{equation}
\end{prop}

For example, let
 $\sigma=5\,\bar{7}\,1\,\bar{4}\,9\,\bar{2}\,\bar{6}\,3\,8$.
 Then we have
$$\rm{Leh}~ \sigma= (1,-2,1,-2,5,-2,-5,3,8)$$
and
$$\expandafter{\rm Rmil_B ~Leh}~\sigma=\expandafter{\rm Rmil_B}~\sigma=\{1,3,8\}, $$
$$\expandafter{\rm Max~Leh}~\sigma=\expandafter{\rm Lmap_B}~\sigma=\{1,5\}.$$

The above proposition implies that the
  $\rm A\text{-}code$ is a bijection from  $B_n$ to ${\rm SE}^\expandafter{\rm B}_n$.
It is easy to see that $\expandafter{\rm{Rmil_B}}~\textbf{i}\sigma=\expandafter{\rm{Lmap_B}}~\sigma$ and $\expandafter{\rm{Rmil_B}}~\sigma=\expandafter{\rm{Lmap_B}}~\textbf{i}\sigma$. So we are led to the  following theorem which shows that the two set-valued statistics $\expandafter{\rm Rmil_B}$ and $\expandafter{\rm Lmap_B}$ for a signed permutation $\sigma$ can be determined by the A-code of $\sigma$.

\begin{theo}\label{A-typeB}
 For any $\sigma\in B_n$, we have
\begin{equation}\label{equA-B}
  \expandafter{(\rm Rmil_B, Lmap_B)}~\sigma=\expandafter{(\rm Max,Rmil_B)~A\text{-}code}~\sigma.
\end{equation}
\end{theo}

Next we define the B-code for a signed permutation. Let $\sigma=\sigma_1\sigma_2\cdots\sigma_n\in B_n$. For $1\leq i \leq n$, let $k_i$ be the smallest integer $k\geq 1$ such that
$|\sigma^{-k}(i)|\leq i$. We define the B-code of $\sigma$ to be the integer sequence $(b_1,b_2,\ldots,b_n)$ with $b_i={(\sigma^{-k_i})}(i)$.
For example, the B-code of the
 signed permutation $\sigma=3\,\bar{1}\,\bar{6}\,\bar{5}\,4\,2$ is
$(1,-1,1,-4,-4,-3)$.

The B-code of a signed
permutation can be also defined recursively as follows. First,
the $\rm
B\text{-}codes$ of the two signed permutations of $B_1$ are defined
as $ \rm B\text{-}code ~1=(1) $ and $\rm B\text{-}code
~\bar{1}=(-1)$. For $n\geq 2$, we write a signed permutation
$\sigma\in B_n$ as a product of disjoint signed
cycles. There are two cases.
\begin{itemize}
\item[Case 1.] Assume that $n$ has a positive sign in $\sigma$ or $\sigma_n=\bar{n}$.
Let $\sigma'\in B_{n-1}$ be the signed permutation obtained from $\sigma$ by deleting $n$ (or $\bar{n}$) in its cycle decomposition.
Here if $n$ (or $\bar{n}$) is in a cycle of length $1$, we just delete this cycle.
Let $b'=(b_1,b_2,\ldots,b_{n-1})$ be the $\rm B\text{-}code$ of $\sigma'$. Then we define the $\rm B\text{-}code$ of $\sigma$ to be $b=(b_1,b_2,\ldots,b_{n-1},\sigma^{-1}(n))$.
\item[Case 2.] Assume that $n$ has a minus sign in $\sigma$ and $\sigma_n\neq\bar{n}$. Changing the sign of $\sigma_n$ and deleting $\bar{n}$ in the cycle decomposition of $\sigma$, we obtain a signed permutation in $B_{n-1}$, denoted by $\sigma'$.
    Let $b'=(b_1,b_2,\ldots,b_{n-1})$ be the $\rm B\text{-}code$ of $\sigma'$. Then we define the $\rm B\text{-}code$ of $\sigma$ to be
    $b=(b_1,b_2,\ldots,b_{n-1},\sigma^{-1}(n))$.
\end{itemize}

The following theorem shows that the set-valued statistics $\rm Lmap_B$
 and $\rm Cyc_B$ of a signed permutation can be computed from the $\textrm{B}$-code.

\begin{theo}\label{B-typeB}
The $\rm B\text{-}code$ is a bijection from $B_n$  to ${\rm SE}^{\rm
B}_n$. Furthermore, for any $\sigma\in B_n$, we have
\begin{equation}\label{equB-B}
  (\expandafter{\rm Cyc_B,Lmap_B)}~\sigma=\expandafter{\rm(Max,Rmil_B)~B\text{-}code}~\sigma.
\end{equation}
\end{theo}

\pf From the recursive definition, it is readily seen that the
$\rm B\text{-}code$ is a bijection from $B_n$ to ${\rm SE}^{\rm
B}_n$. We shall use induction on $n$ to prove (\ref{equB-B}). Clearly, the statement holds for $n=1$.
 Assume that (\ref{equB-B}) holds for $n-1$, where $n\geq 2$. Let $\sigma=\sigma_1\sigma_2\cdots\sigma_n$ be a signed permutation of  $B_n$ with $\expandafter{\rm B\text{-}code}~ b$. Assume that
 $\sigma'$ is the signed permutation of $B_{n-1}$ given in the recursive definition of the B-code. Let $b'=(b_1,b_2,\ldots,b_{n-1})$ be the $\rm B\text{-}code$ of $\sigma'$.

Now we claim that $\expandafter{\rm Cyc_B}~\sigma=\expandafter{\rm Max} ~b$.
There are two cases according to the sign of $n$ in $\sigma$.

First, we consider the case when $n$ has a positive sign in $\sigma$.
       If $\sigma_n\neq n$, let $t=\sigma^{-1}(n)$.
       Since $\sigma'$ is obtained from $\sigma$ by deleting $n$ in its cycle form,
        the $\rm B\text{-}code$ of $\sigma$
        is $b=(b_1,b_2,\ldots,b_{n-1},t)$. Since $0<t<n$,
        we have
        $\expandafter{\rm Cyc_B}~\sigma=\expandafter{\rm Cyc_B}~\sigma' $ and $\expandafter{\rm Max}~b'=\expandafter{\rm Max}~b$. By the induction hypothesis, $\expandafter{\rm Cyc_B}~\sigma'
         =\expandafter{\rm Max}~b'$. Hence ${\rm Cyc_B}~\sigma={\rm Max}~b$.
        If $\sigma_n=n$,
       it can be easily checked that
        $$\expandafter{\rm Cyc_B}~\sigma=\expandafter{\rm Cyc_B}~\sigma'\cup \{n\}=\expandafter{\rm Max}~b'\cup \{n\}=\expandafter{\rm Max}~b. $$

 Then we consider the case when $n$ has a minus sign in $\sigma$.
    If $\sigma_n=\bar{n}$, it is easy to see that
  $$\expandafter{\rm Cyc_B}~\sigma=\expandafter{\rm Cyc_B}~\sigma'
           =\expandafter{\rm Max}~b'=\expandafter{\rm Max}~b. $$
If $\sigma_n\neq\bar{n}$, let $t=\sigma^{-1}(n)$.
Since $n$ has a minus sign in $\sigma$, we have $t<0$.
 Since $b'=(b_1,b_2,\ldots,b_{n-1})$ is the $\rm B\text{-}code$
          of $\sigma'$,  we find that the $\rm B\text{-}code$ of $\sigma$
        is $b=(b_1,b_2,\ldots,b_{n-1},t)$. Since $-n<t<0$, we have
        $\expandafter{\rm Cyc_B}~\sigma=\expandafter{\rm Cyc_B}~\sigma'$ and $\expandafter{\rm Max}~b'=\expandafter{\rm Max}~b$. By
        the induction hypothesis, we get
        $\expandafter{\rm Cyc_B}~\sigma'
         =\expandafter{\rm Max}~b'$. Thus we obtain
        $\expandafter{\rm Cyc_B}~\sigma=\expandafter{\rm Max}~b$.

We now turn to the proof of the relation
$\expandafter{\rm Lmap_B}~\sigma=\expandafter{\rm Rmil_B}~b$. There are four cases.

\noindent Case 1: $\sigma_n=n-1$. By the recursive definition of the B-code, we express $\sigma$ and $\sigma'$ in the one-line notation as follows. For convenience, we display the identity permutation on the top,
        \[
        \begin{array}{ccccccc}
              & 1 & \cdots & |\sigma^{-1}(n)| & \cdots & n-1 & n \\[5pt]
            \sigma= & \sigma_1 & \cdots & v\cdot n & \cdots & \sigma_{n-1} & n-1 \\[5pt]
            \sigma'= &\sigma_1 & \cdots & v\cdot (n-1) & \cdots & \sigma_{n-1}.  &
        \end{array}
        \]
        Here $v=1$ if $n$ has a positive sign in $\sigma$ and $v=-1$ if $n$ has a minus sign in $\sigma$.
        It can be readily seen that $\expandafter{\rm Lmap_B}~\sigma=\expandafter{\rm Lmap_B}~\sigma'$.
        Since $b'=(b_1,b_2,\ldots,b_{n-1})$ is the $\rm B\text{-}code$
          of $\sigma'$, we have $b_{n-1}=\sigma^{-1}(n)$
        and the $\rm B\text{-}code$ of $\sigma$
        is $b=(b_1,b_2,\ldots,b_{n-1},\sigma^{-1}(n))$.
        It follows that
         $\expandafter{\rm Rmil_B}~b=\expandafter{\rm Rmil_B}~b'$.
        By the induction hypothesis, we get
        ${\rm Lmap_B}~\sigma'={\rm Rmil_B}~b' $.
       Hence we deduce that  $\expandafter{\rm Lmap_B}~\sigma=
        \expandafter{\rm Rmil_B}~b$.

\noindent Case 2:  $\sigma_n=\overline{n-1}$.
     If $n$ has a minus sign in $\sigma$, let $t$ be the positive integer such that $\sigma_t=\bar{n}$. As in Case 1, we
     express $\sigma$ and $\sigma'$ as follows
        \[
        \begin{array}{ccccccc}
             & 1 & \cdots & t & \cdots & n-1 & n \\[5pt]
            \sigma= & \sigma_1 & \cdots & \bar{n} & \cdots & \sigma_{n-1} & \overline{n-1} \\[5pt]
            \sigma'= &\sigma_1 & \cdots & n-1 & \cdots & \sigma_{n-1}.  &
        \end{array}
        \]
Clearly, $\expandafter{\rm Lmap_B}~\sigma=\expandafter{\rm Lmap_B}~\sigma'\setminus\{t\}$.
        Since $b'=(b_1,b_2,\ldots,b_{n-1})$ is the $\rm B\text{-}code$
          of $\sigma'$, we have $b_{n-1}={\sigma'}^{-1}(n-1)=t$.
        From the recursive construction of the B-code, it follows that the $\rm B\text{-}code$ of $\sigma$
        is $b=(b_1,b_2,\ldots,b_{n-1},-t\,)$. This implies that
         $\expandafter{\rm Rmil_B}~b=\expandafter{\rm Rmil_B}~b'\setminus\{t\}$.
        By the induction hypothesis, we obtain ${\rm Lmap_B}~\sigma'={\rm Rmil_B}~b' $.
        Therefore $\expandafter{\rm Lmap_B}~\sigma=
        \expandafter{\rm Rmil_B}~ b$.
         If
         $n$ has a positive sign in $\sigma$, let $t$ be the positive integer such that $\sigma_t=n$. Then $\sigma$ and $\sigma'$ can be expressed as follows
        \[
        \begin{array}{ccccccc}
             & 1 & \cdots & t & \cdots & n-1 & n \\[5pt]
            \sigma= & \sigma_1 & \cdots & n & \cdots & \sigma_{n-1} & \overline{n-1} \\[5pt]
            \sigma'= &\sigma_1 & \cdots &\overline{n-1} & \cdots & \sigma_{n-1}.  &
        \end{array}
        \]
      In this case,
        we have $\expandafter{\rm Lmap_B}~\sigma=\expandafter{\rm Lmap_B}~\sigma'\cup\{t\}$.
        Since $b'=(b_1,b_2,\ldots,b_{n-1})$ is the $\rm B\text{-}code$
          of $\sigma'$, then $b_{n-1}=-t$
        and the $\rm B\text{-}code$ of $\sigma$
        is $b=(b_1,b_2,\ldots,b_{n-1},t\,)$.
        It follows that $\expandafter{\rm Rmil_B}~b=\expandafter{\rm Rmil_B}~b'\cup\{t\}$. By the induction hypothesis,
        we deduce that ${\rm Lmap_B}~\sigma'={\rm Rmil_B}~b' $.
        So we arrive at  $\expandafter{\rm Lmap_B}~\sigma=
        \expandafter{\rm Rmil_B}~b$.

   \noindent Case 3:  $\sigma_n\neq n-1$, $\sigma_n\neq \overline{n-1}$ and $|\sigma^{-1}({n-1})|<|\sigma^{-1}(n)|$.
    If $n$ has a positive sign in $\sigma$, let $\sigma_t=n$. Following the similar argument as in Case 2, we have $\expandafter{\rm Lmap_B}~\sigma=\expandafter{\rm Lmap_B}~\sigma'\cup\{t\}$
        and $\expandafter{\rm Rmil_B}~b=\expandafter{\rm Rmil_B}~b'\cup\{t\}$.
        By the induction hypothesis, we deduce that ${\rm Lmap_B}~\sigma'={\rm Rmil_B}~b' $. Hence $\expandafter{\rm Lmap_B}~\sigma=
        \expandafter{\rm Rmil_B}~b$.
    If $n$ has a minus sign in $\sigma$, it can be verified that $\expandafter{\rm Lmap_B}~\sigma=\expandafter{\rm Lmap_B}~\sigma'$
        and $\expandafter{\rm Rmil_B}~b=\expandafter{\rm Rmil_B}~b'$.
        Therefore, we obtain $\expandafter{\rm Lmap_B}~\sigma=
        \expandafter{\rm Rmil_B}~b$.

   \noindent Case 4:  $\sigma_n\neq n-1$, $\sigma_n\neq \overline{n-1}$ and $|\sigma^{-1}({n-1})|>|\sigma^{-1}(n)|$.
    If $n$ has a positive sign in $\sigma$, let $\sigma_t=n$. We write $\sigma$ and $\sigma'$ as follows
    \[
        \begin{array}{ccccccccc}
             & 1   & \cdots & t & \cdots  &|\sigma^{-1}(n-1)| & \cdots & n-1 & n \\[5pt]
            \sigma= & \sigma_1 & \cdots & n & \cdots & v_{n-1}\cdot(n-1) & \cdots &\sigma_{n-1}&\sigma_n \\[5pt]
    \sigma'= & \sigma_1 & \cdots &\sigma_n & \cdots &v_{n-1}\cdot(n-1)  &\cdots &\sigma_{n-1},
        \end{array}
        \]
        where $v_{n-1}=1$ if $n-1$ appears as an element in $\sigma$ and $v_{n-1}=-1$ if $\overline{n-1}$ appears as an element in $\sigma$.
       It can be seen that
        $$\expandafter{\rm Lmap_B}~\sigma=
          (\expandafter{\rm Lmap_B}~\sigma'\cap[1,t-1])\cup\{t\}.$$
        Since $b'=(b_1,b_2,\ldots,b_{n-1})$ is the $\rm B\text{-}code$
          of $\sigma'$, we have $b_{n-1}=\sigma^{-1}(n-1)$
        and the $\rm B\text{-}code$ of $\sigma$
        is $b=(b_1,b_2,\ldots,b_{n-1},t)$.
        Hence we get $$\expandafter{\rm Rmil_B}~b=
        (\expandafter{\rm Rmil_B}~b'\cap[1,t-1])\cup\{t\}.$$
        By the induction hypothesis, we obtain ${\rm Lmap_B}~\sigma'={\rm Rmil_B}~b' $.  Thus we get $\expandafter{\rm Lmap_B}~\sigma=
        \expandafter{\rm Rmil_B}~b$.
     If $n$ has a minus sign in $\sigma$, it can be checked that
         $\expandafter{\rm Lmap_B}~\sigma=
          \expandafter{\rm Lmap_B}~\sigma'\cap[1,-\sigma^{-1}(n)-1]$
        and $\expandafter{\rm Rmil_B}~b=
        \expandafter{\rm Rmil_B}~b'\cap[1,-\sigma^{-1}(n)-1]$.
        By the induction hypothesis, we conclude that $\expandafter{\rm Lmap_B}~\sigma=
        \expandafter{\rm Rmil_B}~b$. This completes the proof. \qed

In fact, it can be shown that the pair of statistics $(inv_B,nmin_B)$  of a signed permutation $\sigma$ can be recovered from its A-code and the pair of statistics $(sor_B,l_B')$
can be recovered from its B-code.

We now describe how to
recover a signed permutation $\sigma$ from its A-code $a=(a_1,a_2,\ldots,a_n)\in
\textrm{SE}^\textrm{B}_n$. It is essentially the
same as the procedure to
recover a permutation from the inversion code.

We start with the empty word $\sigma^{(0)}$,
 then it will take $n$ steps
to construct a signed permutation $\sigma$ with A-code $a$.
At the first step, if $a_1=1$, then set $\sigma^{(1)}=1$. If $a_1=-1$, then set $\sigma^{(1)}=\bar{1}$.
For $1< i\leq n$, assume that at step $i$, we have constructed a signed permutation
$\sigma^{(i-1)}\in B_{i-1}$. If $|a_i|=1$, the signed permutation $\sigma^{(i)}$ is obtained  by inserting
the element  $i$ with a sign of $a_i$ before the first element
of $\sigma^{(i-1)}$.
If $|a_i|>1$,
then the signed permutation $\sigma^{(i)}$ is obtained from $\sigma^{(i-1)}$ by inserting
the element  $i$ with a sign of $a_i$ immediately after the
 $(|a_i|-1)$-th element in $\sigma^{(i-1)}$.
Eventually, the signed permutation $\sigma^{(n)}$ is  a signed permutation $\sigma$  with A-code $a$.
For example, $a=(1,1,-3,-2,3)$, then we have
\[
\begin{array}{lll}
          & \sigma^{(0)}=  & \emptyset, \\[5pt]
  a_1=1,  & \sigma^{(1)}=  & 1, \\[5pt]
  a_2=1,  & \sigma^{(2)}=  & 2~1, \\[5pt]
  a_3=-3, & \sigma^{(3)}=  & 2~1~\bar{3}, \\[5pt]
  a_4=-2, & \sigma^{(4)}=  & 2~\bar{4}~1~\bar{3}, \\[5pt]
  a_5=3,  & \sigma^{(5)}=  & 2~\bar{4}~5~1~\bar{3}.
\end{array}
\]
So the signed permutation $2\,\bar{4}\,5\,1\,\bar{3}$ corresponds to
the A-code $(1,1,-3,-2,3)$.

The relationship between a signed permutation $\sigma$ and its B-code $b=(b_1,b_2,\ldots,b_n)$ can be described as follows. Let $\sigma'$ be the signed permutation obtained from $\sigma$ as in the recursive construction of the B-code.
So the B-code of $\sigma'$ is $b'=(b_1,b_2,\ldots,b_{n-1})$.
If $n$ has a positive sign in $\sigma$ or $\sigma_n=\bar{n}$, then $\sigma'$ is obtained from $\sigma$ by deleting $n$ in its cycle decomposition. Let $(i,i)$ denote the identity permutation for any $1\leq i\leq n$. Since $b_n=\sigma^{-1}(n)$, we have $\sigma=\sigma'(b_n,n)$. We note here that $\sigma'$ is considered as a signed permutation of $B_n$ which maps $n$ to $n$.
If $n$ has a minus sign in $\sigma$ and $\sigma_n\neq \bar{n}$, then
 $\sigma'$ is obtained from $\sigma$ by changing the sign of $\sigma_n$ and deleting $\bar{n}$ in its cycle decomposition.
 Since $b_n=\sigma^{-1}(n)$,
 it is readily seen that  $\sigma=\sigma'(b_n,n)$. Again here $\sigma'$ is considered as a signed permutation of $B_n$ which maps $n$ to $n$.
 Hence we obtain that
 $\sigma=(b_1,1)(b_2,2)\cdots(b_n,n)$.

The following lemma gives  expressions of $inv_B(\sigma)$
and $nmin_B(\sigma)$ in terms of the A-code of $\sigma$.

\begin{lemm}\label{lem-a-typeB}
 For a signed permutation $\sigma\in B_n$ with {\rm A-code} $a=(a_1,a_2,\ldots,a_n)$, we have
    \begin{equation}\label{inv-b}
      inv_B(\sigma)=\sum^n_{i=1}(i-a_i-\chi(a_i<0))
    \end{equation}
   and
    \begin{equation}\label{nmin-b}
       ~ nmin_B(\sigma)=n-|{\rm Max}~a|.
    \end{equation}
\end{lemm}

\pf
Consider the procedure to recover a signed permutation
 from the A-code $a$. It is easily seen that after the
 $i$-th step, the type $B$ inversion number  increases by
  $i-a_i$ when $a_i>0$  and by $i-a_i-1$ when $a_i<0$.
Hence we have
$$inv_B(\sigma^{(i)})-inv_B(\sigma^{(i-1)})=i-a_i-\chi(a_i<0). $$
 Since $inv_B(\sigma^{(0)})=0$, we find $$inv_B(\sigma)=\sum^n_{i=1}(i-a_i-\chi(a_i<0)).$$
 In view of (\ref{equA-B}), it is easy to see that $nmin_B(\sigma)=n-rl\textrm{-}min_B(\sigma)=n-|\textrm{Rmil}_{\rm B}~\sigma|=n-|\textrm{Max}~a|$. This completes the proof. \qed

The following lemma shows that $sor_B(\sigma)$
and $l_B'(\sigma)$ can be expressed in terms of the B-code of $\sigma$.

\begin{lemm}\label{lem-b-typeB}
  For a signed permutation $\sigma\in B_n$ with {\rm B-code} $b=(b_1,b_2,\ldots,b_n)$, we have
  \begin{equation}\label{sor-b}
    sor_B(\sigma)=\sum^n_{i=1}(i-b_i-\chi(b_i<0))
  \end{equation}
 and
  \begin{equation}\label{l-b}
    l_B'(\sigma)=n-|{\rm Max}~b|.
  \end{equation}
\end{lemm}

\pf
Since $b=(b_1,b_2,\ldots,b_n)$ is the B-code of $\sigma$, it has
been shown that
$$\sigma=(b_1,1)(b_2,2)\cdots(b_n,n).$$
 By the definition of the sorting index of $\sigma$, we  see that
$$sor_B(\sigma)=\sum^n_{i=1}(i-b_i-\chi(b_i<0)).$$
From (\ref{equB-B}) it follows that
 $l_B'(\sigma)=n-cyc_B(\sigma)=n-|{\rm Cyc}_{\rm B}~\sigma|=n-|\textrm{Max}~b|$.
 This completes the proof. \qed

Combining Theorem \ref{A-typeB}, Theorem \ref{B-typeB}, Lemma \ref{lem-a-typeB} and Lemma \ref{lem-b-typeB}, we obtain  the equidistribution of $(inv_B,{\rm Lmap_B},{\rm Rmil_B})$ and
$(sor_B,{\rm Lmap_B},{\rm Cyc_B})$ over $B_n$.

\begin{theo}\label{thm-set-B}
The map $\psi\colon B_n\longrightarrow B_n$ defined by
$\psi=\rm (B\textrm{-}code)^{-1}\circ A\textrm{-}code$ is a bijection. For any $\sigma\in B_n$, we have
\begin{equation}
(inv_B,{\rm Lmap_B},{\rm Rmil_B})~\sigma=(sor_B,{\rm Lmap_B},{\rm Cyc_B})~\psi(\sigma).
\end{equation}
In particular,
\begin{equation}
(inv_B,nmin_B)~\sigma=(sor_B,l_B')~\psi(\sigma).
\end{equation}
\end{theo}

Notice that
$\expandafter{\rm Cyc_B}~\sigma=\expandafter{\rm Cyc_B}~\textbf{i}\sigma $ and
$\expandafter{\rm Lmap_B}~\sigma=\expandafter{\rm Rmil_B}~\textbf{i}\sigma$.
Thus Theorem \ref{thm-set-B} implies the
 following equidistributions  which can be
 viewed as type $B$ analogues of the equidistributions
  given in Theorem \ref{thm-set-A}.

\begin{theo}
The six pairs of set-valued statistics ~$\rm (Cyc_B,Rmil_B)$,
$\rm(Cyc_B,Lmap_B)$, $\rm(Rmil_B,Lmap_B)$, $\rm(Lmap_B,Rmil_B)$,
$\rm(Lmap_B,Cyc_B)$ and $\rm(Rmil_B,Cyc_B)$ are equidistributed over
$B_n$:
\end{theo}
 \begin{center}
 \begin{tabular}{ccccccccccc}
 $B_n$ &  $\stackrel{ \textbf{i}}{\longrightarrow}$ &  $B_n$   &  $\stackrel{\psi^{-1}}{\longrightarrow}$    & $ B_n $ & $ \stackrel{\textbf{i}}{\longrightarrow}$ & $ B_n $ & $\stackrel{\psi }{\longrightarrow}$ & $ B_n $ &
$\stackrel{\textbf{i}}{\longrightarrow}$ &$ B_n $ \\[5pt]
    $\rm Cyc_B\choose Rmil_B $   &         & $\rm Cyc_B \choose Lmap_B $   &   & $\rm Rmil_B\choose Lmap_B $  &      & $\rm Lmap_B\choose Rmil_B  $ & & $\rm Lmap_B\choose Cyc_B$ &    &$\rm Rmil_B\choose Cyc_B$.
 \end{tabular}
\end{center}

The above theorem for set-valued statistics reduces to the following
equidistributions of pairs of statistics of signed permutations.
It is clear that  $nmin_B(\sigma)=
nmax_B(\textbf{i}\sigma)$. Since the bijection $\psi$ preserves $\textrm{Lmap}_\textrm{B}$, it is easy to see
that $\psi$ also preserves the statistic $nmax_B$.
Hence we are led to the following assertion.

\begin{coro}
The four pairs of statistics $(sor_B,l'_B)$, $(inv_B,nmin_B)$,
$(inv_B,nmax_B)$ and $(sor_B,nmax_B)$ are
equidistributed over $B_n$:
\end{coro}
\begin{center}
 \begin{tabular}{ccccccc}
 $B_n$ &  $\stackrel{\psi^{-1}}{\longrightarrow}$ &  $B_n$   &  $\stackrel{\textbf{i} }{\longrightarrow}$    & $ B_n $ & $ \stackrel{\psi}{\longrightarrow}$ & $ B_n $\\[5pt]
   $sor_B\choose l'_B $   &         & $ inv_B\choose nmin_B $   &   & $inv_B\choose nmax_B $  &      & $sor_B\choose nmax_B  $.
 \end{tabular}
\end{center}

\section{A bijection on $D_n$}

In this section, we define two statistics $nmin_D$ and $\tilde{l}_D'$ for elements of a Coxeter group of type $D$ and we construct a bijection to derive the equidistribution of  the pairs of statistics $(inv_D,nmin_D)$ and $(sor_D,\tilde{l}_D')$. This yields a refinement of Petersen's equidistribution of $inv_D$ and $sor_D$.

The type $D$ Coxeter group $D_n$ is the subgroup of $B_n$ consisting
of  signed permutations with an even number of minus signs in the
signed permutation notation. As a set of generators for $D_n$, we take
$$ S^D=\{(\bar{1},2),(1,2),(2,3),\ldots,(n-1,n)\}. $$
For simplicity, let $s_i=(i, i+1)$ for $1\leq i<n$ and $s_{\bar{1}}=(\bar{1},2) $. The set of
reflections of $D_n$ is
$$R^D=\{(i,j):1\leq |i|<j\leq n\}. $$
For $\sigma=\sigma_1\sigma_2\cdots\sigma_n\in D_n$, the type $D$ inversion number of $\sigma$ is given by
$$inv_D(\sigma)=|\{(i,j):1\leq i<j\leq n,\sigma_i>\sigma_j\}|+|\{(i,j):1\leq i<j\leq n,\overline{\sigma_i}>\sigma_j\}|.$$
The length of $\sigma$, denoted $l_D(\sigma)$, is
the minimal number of transpositions in $S^D$ needed to express
$\sigma$. It is known that
$l_D(\sigma)=inv_D(\sigma)$, see \cite[Section 8.2]{BB}.

It is well-known that the generating function of $l_D$ is
\begin{equation}
   \sum_{\sigma\in
  D_n}q^{l_D(\sigma)}={[n]}_q\prod_{r=1}^{n-1}{[2r]}_q,
\end{equation}
see \cite{BB}.

Recall that the set of reflections of $B_n$ is $$T^B=\{(i,j):1\leq i<j\leq n\}\cup\{(\bar{i},j):1\leq i\leq j\leq n\}.  $$
For $\sigma\in D_n$,
it has a unique
factorization into  a product of signed transpositions in
$T^B$:
\begin{equation}\label{fac-B}
  \sigma=(i_1,j_1)(i_2,j_2)\cdots(i_k,j_k),
\end{equation}
where $0<j_1<j_2<\cdots<j_k\leq n$.
Petersen defined the type $D$ sorting index of $\sigma$ as
$$sor_D(\sigma)=\sum_{r=1}^k(j_r-i_r-2\chi(i_r<0)).  $$
It has been shown by Petersen that $sor_D$ has the same generating function as $inv_D$.

 \begin{theo}\label{th4.1}
 For $n\geq 4$,
 \begin{equation}
  \sum_{\sigma\in
  D_n}q^{sor_D(\sigma)}={[n]}_q\prod_{r=1}^{n-1}{[2r]}_q.
   \end{equation}
  Thus, $sor_D$ is Mahonian.
  \end{theo}

Next we define two statistics $\tilde{l}_D'$ and $nmin_D$ for a signed permutation $\sigma \in D_n$.
For $1\leq |i|<j\leq n$, we adopt the notation
$t_{ij}$ for the transposition $(i,j)$. For $1<i\leq n$, we define $t_{\bar{i}i}=(\bar{i},i)(\bar{1},1)$. Then we set
$$T^D=\{t_{ij}:1\leq |i|<j\leq n\}\cup\{t_{\bar{i}i}:1<i\leq n \}.$$
 We denote
 by $\tilde{l}_D'(\sigma)$  the
minimal number of elements in $T^D$ that are needed  to express $\sigma$.
Define the statistic $nmin_D$ as follows
$$nmin_D(\sigma)=|\{i:\sigma_i>|\sigma_j|
~\text{for some}~ j>i\}|+N(\sigma\backslash\{\bar{1}\}),$$
 where $N(\sigma\backslash\{\bar{1}\})$ is the number of minus signs associated with elements  greater than  $1$ in
 the signed permutation notation of $\sigma$.

The following theorem is a refinement of the
equidistribution of $inv_D$ and $sor_D$. We shall give a combinatorial proof and an algebraic proof.

\begin{theo}\label{th4.2}
For $n\geq 2$, the two pairs of statistics $(inv_D,nmin_D)$ and $(sor_D,\tilde{l}_D')$ are
equidistributed over $D_n$. Moreover, we have
\begin{align}
\label{eqinvd}
\sum_{\sigma\in
D_n}q^{inv_D(\sigma)}t^{nmin_D(\sigma)}=\prod_{r=1}^{n-1}(1+q^rt+qt\cdot{[2r]}_q)
  ,\\
\label{eqsord}
\sum_{\sigma\in
D_n}q^{sor_D(\sigma)}t^{\tilde{l}_D'(\sigma)}=\prod_{r=1}^{n-1}(1+q^rt+qt\cdot{[2r]}_q).
\end{align}
\end{theo}

To give a combinatorial proof of the equidistribution of $(inv_D,nmin_D)$ and $(sor_D,\tilde{l}_D')$ in Theorem \ref{th4.2}, we
introduce the co-sorting index $sor_D'$ which turns out
to be  equivalent to the sorting index $sor_D$.
To define the co-sorting index,
we need the factorization of an element $\sigma\in D_n$
into elements in $T^D$. More precisely, similarly as (\ref{fac-sigma-B}), we can uniquely express $\sigma\in D_n$ as
$$\sigma=t_{i_1j_1}t_{i_2j_2}~\cdots~t_{i_mj_m}, $$
where $1<j_1<j_2<\cdots<j_m\leq n$. For example, let
$\sigma=\bar{2}\,\bar{4}\,5\,\bar{1}\,\bar{3}$.
Then we have $\sigma=t_{12}t_{\bar{3}3}t_{\bar{2}4}t_{35}$.
Then the co-sorting index of $\sigma$ is defined by
 $$sor_D'(\sigma) =\sum_{r=1}^m(j_r-i_r-2\chi(i_r<0)).$$

\begin{lemm}\label{sor-D}
  For any $\sigma\in D_n$, we have
  $sor_D(\sigma)=sor_D'(\sigma).$
\end{lemm}

\pf
Write $\sigma$ in the following form
\begin{equation}\label{fac}
\sigma=t_{i_1j_1}t_{i_2j_2}~\cdots~t_{i_mj_m},
\end{equation}
where $t_{i_1j_1}, t_{i_2j_2}, \ldots, t_{i_mj_m}\in T^D$ and $1<j_1<j_2<\cdots<j_m\leq n$.
Since the co-sorting index of $\sigma$ can be expressed in terms of the factorization (\ref{fac}), to prove the
 the equivalence of the sorting index  and the co-sorting index of $\sigma$, we wish to rewrite (\ref{fac}) as a product of transpositions in $T^B$ from which
   the sorting index of $\sigma$  can be determined.

In fact, it can be shown that
$\sigma$ can be written as a product of
 transpositions in $T^B$ which is either of the form
\begin{equation}\label{fac-sigma1}
  (p_1,j_1)(p_2,j_2)\cdots(p_m,j_m),
\end{equation}
or of the form
\begin{equation}\label{fac-sigma2}
  (\bar{1},1)(p_1,j_1)(p_2,j_2)\cdots(p_m,j_m),
\end{equation}
where for $1\leq k\leq m$,
\begin{equation}\label{ex-p}
 p_k=\left\{
\begin{array}{llll}
1 ~\text{or}~ \bar{1},& \text{if}~i_k=1,\\[5pt]
1 ~\text{or}~ \bar{1},& \text{if}~i_k=\bar{1},\\[5pt]
i_k, & \text{otherwise}.
\end{array}
\right.
\end{equation}
 To this end, we claim that for $1\leq r\leq m$,
$t_{i_rj_r}t_{i_{r+1}j_{r+1}}\cdots t_{i_mj_m}$ can be expressed as
a product of transpositions in $T^B$ which is either of the form \begin{equation}\label{form1}
  (p_r,j_r)(p_{r+1},j_{r+1})\cdots(p_m,j_m)
\end{equation}
 or of the form
\begin{equation}\label{form2}
  (\bar{1},1)(p_r,j_r)(p_{r+1},j_{r+1})\cdots(p_m,j_m),
\end{equation}
where $p_k$ is given as in (\ref{ex-p}).
Let us first consider the case $r=m$.
In this case,
if $i_m\neq \overline{j_m}$, then
$t_{i_mj_m}$ equals $(i_m,j_m)$ which is of the form
(\ref{form1}).
If $i_m=\overline{j_m}$, then $t_{i_mj_m}$ equals $(\bar{1},1)(i_m,j_m)$ which
is of the form (\ref{form2}).

 Assume the claim holds for $r$, where $1<r\leq m$. We aim to show that it holds for $r-1$.
If $t_{i_rj_r}t_{i_{r+1}j_{r+1}}\cdots t_{i_mj_m}$ can be expressed
 in the form (\ref{form1}),
then we have
$$t_{i_{r-1}j_{r-1}}t_{i_rj_r}\cdots t_{i_mj_m}=\left\{
        \begin{array}{llll}
        (\bar{1},1)(i_{r-1},j_{r-1})(p_r,j_r)\cdots(p_m,j_m),     & \text{if}~i_{r-1}=\overline{j_{r-1}},\\[5pt]
          (i_{r-1},j_{r-1})(p_r,j_r)\cdots(p_m,j_m), & \text{otherwise,}
        \end{array}
      \right.
$$
which is either of the form (\ref{form2}) or of the form (\ref{form1}).
We now assume  that  $t_{i_rj_r}t_{i_{r+1}j_{r+1}}\cdots t_{i_mj_m}$ can be
 expressed in the form (\ref{form2}).
It follows that
 $$ t_{i_{r-1}j_{r-1}}t_{i_rj_r}\cdots t_{i_mj_m}=\left\{
        \begin{array}{llll}
              (i_{r-1},j_{r-1})(p_r,j_r)\cdots(p_m,j_m),      & \text{if}~i_{r-1}=\overline{j_{r-1}},\\[5pt]
          (\bar{1},1)(\overline{i_{r-1}},j_{r-1})(p_r,j_r)\cdots(p_m,j_m), &           \text{if}~i_{r-1}=1~ \text{or}~ \bar{1},\\[5pt]
          (\bar{1},1)(i_{r-1},j_{r-1})(p_r,j_r)\cdots(p_m,j_m), & \text{otherwise,}
        \end{array}
      \right.
$$
which is either of the form (\ref{form1}) or of the form (\ref{form2}). Thus we have verified that the claim holds for
any $1\leq  r\leq m$.

Now we have shown that $\sigma$ can
be expressed as (\ref{fac-sigma1}) or (\ref{fac-sigma2}).
 Then the sorting index
$sor_D(\sigma)$ can be determined by this factorization, namely,
$$sor_D(\sigma)=\sum_{r=1}^m(j_r-p_r-2\chi(p_r<0)).$$
 By (\ref{ex-p}), we find that $$j_r-p_r-2\chi(p_r<0)=j_r-i_r-2\chi(i_r<0)$$
for $1\leq r\leq m$. In view of (\ref{fac}), we see that
$$sor_D'(\sigma)=\sum_{r=1}^m(j_r-i_r-2\chi(i_r<0)).$$
It follows that $sor_D(\sigma)=sor_D'(\sigma)$. This completes the proof. \qed

To justify the equidistribution of $(inv_D,nmin_D)$ and $(sor_D,\tilde{l}'_D)$,
we shall give a bijection which transforms $(inv_D,nmin_D)$ to $(sor_D,\tilde{l}'_D)$.
The bijection will be described in terms of two codes,
called the E-code and the F-code of an element of $D_n$.
 It will be shown that  the pair of statistics ($inv_D,nmin_D)$ can be computed from the E-code  whereas the pair of statistics $(sor_D,\tilde{l}_D')$  can be computed from the F-code.

 Given an element $\sigma\in D_n$, the E-code of $\sigma$ is an integer sequence $e=(e_1,e_2,\ldots,e_n)$ generated
 by the following procedure.
   We wish to construct a sequence of signed permutations $\sigma^{(n)},\sigma^{(n-1)},\ldots,\sigma^{(1)}$ where
     $\sigma^{(i)}\in D_i$ for $1\leq i\leq n$.
   First, we set $\sigma^{(n)}=\sigma$.
 For $i$ from $n$ to $2$, we construct $\sigma^{(i-1)}$ from
 $\sigma^{(i)}$. Consider the letter $i$ in $\sigma^{(i)}$.  If $i$ has a positive sign in $\sigma^{(i)}$, then assume that
  $i$ appears at the $p$-th position in $\sigma^{(i)}$.
  In this case, we set $e_i=p$ and let  $\sigma^{(i-1)}$
   be the signed permutation obtained from $\sigma^{(i)}$ by
  deleting the element $i$.
    If $i$ has a minus sign in $\sigma^{(i)}$, then
we assume that $\bar{i}$ appears at the $p$-th position in
  $\sigma^{(i)}$. We then set $e_i=-p$. Let
   $\sigma'$ be the signed permutation obtained from $\sigma^{(i)}$ by deleting $\bar{i}$, and let $\sigma^{(i-1)}$ be the signed permutation obtained from $\sigma'$ by changing the sign of the element at the first position.
 It can be seen that the resulting signed permutation $\sigma^{(1)}$ is the identity permutation $1$.
 Finally, we set $e_1=1$.
 For example, let $\sigma=2\,\bar{4}\,5\,1\,\bar{3}$. Then we have
\[
\begin{array}{cll}
  \sigma^{(5)}= & 2~\bar{4}~\textbf{5}~1~\bar{3}, & e_5=3, \\[5pt]
  \sigma^{(4)}= & 2~\bar{\textbf{4}}~1~\bar{3}, & e_4=-2, \\[5pt]
  \sigma^{(3)}= & \bar{2}~1~\bar{\textbf{3}}, & e_3=-3, \\[5pt]
  \sigma^{(2)}= & \textbf{2}~1, & e_2=1, \\[5pt]
  \sigma^{(1)}= & 1, & e_1=1. \\
\end{array}
\]
Hence the E-code of $\sigma=2\,\bar{4}\,5\,1\,\bar{3}$ is $(1,1,-3,-2,3)$.

It can be checked that
the above procedure is reversible. In other words,
   one can recover an element $\sigma\in D_n$ from an E-code $e=(e_1,e_2,\ldots,e_n)$.
For $1<r\leq n$, it is routine to verify that \begin{equation}
  inv_D(\sigma^{(r)})-inv_D(\sigma^{(r-1)})=r-e_r-2\chi(e_r<0)
\end{equation} and
\begin{equation}
  nmin_D(\sigma^{(r)})-nmin_D(\sigma^{(r-1)})=1-\chi(e_r=r).
\end{equation}
So we are led to the following formulas for $inv_D(\sigma)$ and
$nmin_D(\sigma)$.

\begin{prop}\label{prop4.3}
Given an element $\sigma\in D_n$, let $e=(e_1,e_2,\ldots,e_n)$ be its {\rm E-code}. Then we have
\begin{equation}
  inv_D(\sigma)=\sum_{r=1}^n (r-e_r-2\chi(e_r<0))
\end{equation}
and
\begin{equation}
  nmin_D(\sigma)=n-\sum_{r=1}^n \chi(e_r=r).
\end{equation}
\end{prop}

We now define the F-code of an element $\sigma\in D_n$ as  an integer sequence $f=(f_1,f_2,\ldots,f_n)$ determined by the following procedure. To compute the F-code $f$, we will generate a sequence of signed permutations $\sigma^{(n)},\sigma^{(n-1)},\ldots,\sigma^{(1)}\in D_n$.  Let us begin with $\sigma^{(n)}=\sigma$.
 For $i$ from $n$ to $2$, we construct
  $\sigma^{(i-1)}$ from
 $\sigma^{(i)}$. Consider the letter $i$ in $\sigma^{(i)}$.
 If $i$ has a positive sign in $\sigma^{(i)}$, say $\sigma^{(i)}(p)=i$, then let $f_{i}=p$ and let  $\sigma^{(i-1)}$ be the signed permutation obtained from  $\sigma^{(i)}$ by exchanging the letter $i$ and the letter at the $i$-th position.
 If $i$ has a minus sign in $\sigma^{(i)}$ and
  $\sigma^{(i)}(i)=\bar{i}$, then let $f_{i}=-i$ and let
  $\sigma^{(i-1)}$ be the signed permutation obtained from
  $\sigma^{(i)}$ by changing both the signs of the element
  at the $i$-th position and the element at the first position.
  If $i$ has a minus sign in $\sigma^{(i)}$
   and $\sigma^{(i)}(i)\neq\overline{i}$,
  say $\sigma^{(i)}(p)=\overline{i}$, then  let $f_{i}=-p$ and let
     $\sigma^{(i-1)}=\sigma^{(i)}(\bar{p},i)$.
   It can be readily seen that the resulting signed permutation $\sigma^{(1)}$
   is the identity permutation $1\,2\cdots n$. Finally, we set $f_1=1$.
  For example, let
$\sigma=\bar{2}\,\bar{4}\,5\,\bar{1}\,\bar{3}$. Then we proceed as follows
\[
\begin{array}{cclc}
  \sigma^{(5)}= & \bar{2}~\bar{4}~\textbf{5}~\bar{1}~\bar{3},& & f_5=3, \\[5pt]
  \sigma^{(4)}= & \bar{2}~\bar{\textbf{4}}~\bar{3}~\bar{1}~5,& & f_4=\bar{2}, \\[5pt]
  \sigma^{(3)}= & \bar{2}~1~\bar{\textbf{3}}~4~5,& & f_3=\bar{3}, \\[5pt]
  \sigma^{(2)}= & \textbf{2}~1~3~4~5,& & f_2=1, \\[5pt]
  \sigma^{(1)}= & 1~2~3~4~5,& & f_1=1. \\
\end{array}
\]
Hence the F-code of $\sigma=\bar{2}\,\bar{4}\,5\,\bar{1}\,\bar{3}$ is $(1,1,-3,-2,3)$.
It is easily seen that the above procedure is reversible. So we can
recover $\sigma$ from its F-code.

The following proposition gives  expressions  of $sor_D$ and $\tilde{l}'_D(\sigma)$ in terms of the F-code.

\begin{prop}\label{prop4.4}
Given an element $\sigma\in D_n$, let $f=(f_1,f_2,\ldots,f_n)$ be its {\rm F-code}. Then we have
\begin{equation}\label{sor-D-expression}
  sor_D(\sigma)=\sum_{r=1}^n (r-f_r-2\chi(f_r<0))
\end{equation}
and
\begin{equation}\label{MRL-D}
  \tilde{l}'_D(\sigma)=n-\sum_{r=1}^n \chi(f_r=r).
\end{equation}
\end{prop}

\pf
For $1\leq i\leq n$, we let $t_{ii}$ denote the identity permutation.
Examining  the procedure to construct the F-code
of $\sigma$, we see that for $1<r\leq n$, we have
\begin{equation}\label{relation-R}
  \sigma^{(r)}=\sigma^{(r-1)} t_{f_rr}.
\end{equation}
 It follows that
\begin{equation}\label{fac-D-R1}
  \sigma^{(r)}=t_{f_11}
t_{f_22}~\cdots~t_{f_rr}.
\end{equation}
By the definition of the co-sorting index, we find
\begin{equation}
  sor_D'(\sigma^{(r)})-sor_D'(\sigma^{(r-1)})=r-f_r-2\chi(f_r<0).
\end{equation}
Applying Lemma \ref{sor-D}, we get
\begin{equation} \label{sd}
  sor_D(\sigma^{(r)})-sor_D(\sigma^{(r-1)})=r-f_r-2\chi(f_r<0).
\end{equation}
Summing (\ref{sd}) over $r$ gives (\ref{sor-D-expression}).

To prove (\ref{MRL-D}),
it suffices to show that
\begin{equation}\label{MRL-R}
  \tilde{l}'_D(\sigma^{(r)})-\tilde{l}'_D(\sigma^{(r-1)})=1-\chi(f_r=r)
\end{equation}
for $1<r\leq n$.
If $f_r=r$, it is clear that $\sigma^{(r)}=\sigma^{(r-1)}$. So (\ref{MRL-R}) holds in this case.
If $f_r\neq r$, let  $\tilde{l}'_D(\sigma^{(r)})=l$.
 Then $\sigma^{(r)}$ can be decomposed as
 \begin{equation}\label{fac-D-R}
  \sigma^{(r)}=t_{i_1j_1}t_{i_2j_2}~\cdots~t_{i_lj_l}
\end{equation}
where $t_{i_1j_1},t_{i_2j_2},\cdots,t_{i_lj_l}\in T^D$.
For $t=t_{ij}\in T^D$ and $1<k\leq n$, we say that
$t$ fixes $k$
 if and only if $k\neq i,\bar{i},j ~\text{or}~\bar{j}$ in the sense that if $k\neq i,\bar{i},j ~\text{or}~\bar{j}$, then $t_{ij}$ maps $k$ to $k$ when we consider $t_{ij}$ as a map on $\{1,2,\ldots,n,\bar{1},\bar{2},\ldots,\bar{n}\}$.
It can be verified that for any $1<k\leq n$ and $t_1,t_2\in T^D$, there exist $t_3,t_4\in T^D$ such that $t_1t_2=t_3t_4$ and
$t_3$ fixes $k$.
Thus we can use (\ref{fac-D-R}) to derive an expression of $\sigma^{(r)}$ of the form \begin{equation}\label{fac-D-R2}
  \sigma^{(r)}=t_{i_1'j_1'}t_{i_2'j_2'}~\cdots~t_{i_l'j_l'}
\end{equation}
 where $t_{i_1'j_1'},t_{i_2'j_2'},\cdots,t_{i_l'j_l'}\in T^D$ and
$t_{i_p'j_p'}$ fixes $r$ for $1\leq p\leq l-1$.
Since $f_r\neq r$, it follows from (\ref{fac-D-R1}) that $\sigma^{(r)}$ maps $f_r$ to $r$. Hence we deduce  that $t_{i_l'j_l'}=t_{f_rr}$.
By (\ref{relation-R}) and (\ref{fac-D-R2}), we get
$$t_{i_1'j_1'}t_{i_2'j_2'}\cdots t_{i_{l-1}'j_{l-1}'}=\sigma^{(r-1)}.$$ Hence we arrive at $$\tilde{l}'(\sigma^{(r-1)})\leq l-1.$$
From (\ref{relation-R}), it is clear that
$$l\leq \tilde{l}'(\sigma^{(r-1)})+1,$$ so we conclude that
\begin{equation}
  l=\tilde{l}'(\sigma^{(r-1)})+1.
\end{equation}
This completes the  proof of (\ref{MRL-D}).
\qed

Combining Propositions \ref{prop4.3} and \ref{prop4.4}, we obtain
a bijection $\rho\colon D_n\longrightarrow D_n$ which leads to the equidistrubution of $(inv_D,nmin_D)$ and $(sor_D,\tilde{l}_D')$.
More precisely, the bijection $\rho$ is given by
$$\rho={\rm F\text{-}code}^{-1}\circ {\rm E\text{-}code}.$$
Then we arrive at the following theorem.

\begin{theo}\label{th4.5}
 The bijection $\rho$ transforms
 $(inv_D,nmin_D)$ to $(sor_D,\tilde{l}_D')$, that is, for any $\sigma\in D_n$, we have
\begin{equation}\label{equation-D}
  (inv_D,nmin_D)~\sigma=~(sor_D,\tilde{l}'_D)~\rho(\sigma).
\end{equation}
\end{theo}

We now  present a proof of Theorem \ref{th4.2}
based on two different factorizations of the diagonal sum
$\sum_{\sigma\in D_n}\sigma$  in the group algebra $\mathbb{Z}[D_n]$.
It turns out that the bivariate  generating functions of  $(inv_D,nmin_D)$ and $(sor_D,\tilde{l}_D')$ are both equal to
$$D_n(q,t)=\prod_{r=1}^{n-1}(1+q^rt+qt\cdot{[2r]}_q).$$

To derive the bivariate generating function of $(inv_D,nmin_D)$, we recall Petersen's  factorization of the diagonal sum $\sum_{\sigma\in D_n}\sigma$.
The elements
$\Psi_1,\Psi_2,\ldots,\Psi_{n-1}$ of the group algebra of $D_n$ are recursively defined as follows.
Recall that $s_i=(i, i+1)$ for $1\leq i<n$ and $s_{\bar{1}}=(\bar{1},2)$.
For $i=1$, let
$$\Psi_1=1+s_1+s_{\bar{1}}+s_1s_{\bar{1}}.$$
For $i\geq2$, let
$$\Psi_i=1+s_i\Psi_{i-1}+s_i\cdots s_2 s_1 s_{\bar{1}}s_2\cdots s_i.$$
Petersen found the following factorization.

\begin{prop} \label{GA1} For $n\geq 2$,
  $$\sum_{\sigma\in D_n}\sigma=\Psi_1\Psi_2\cdots\Psi_{n-1}. $$
\end{prop}

For an element $\sigma \in D_n$, we define  the weight of $\sigma$ to be
$$\mu(\sigma)=q^{inv_D(\sigma)}t^{nmin_D(\sigma)}.$$ As usual, the weight
function is considered as a linear map on $\mathbb{Z}[D_n]$. It can be routinely checked that
\begin{equation}\label{equ-psi}
  \mu(\Psi_i)=1+tq^i+tq\,(1+q+\cdots+q^{2i-1})=1+tq^i+tq\,[2i]_q.
\end{equation}

We are now ready to finish the proof of
 relation (\ref{eqinvd}) concerning the  bivariate generating function of $(inv_D,nmin_D)$.

\noindent
{\em Proof of (\ref{eqinvd}) in Theorem \ref{th4.2}.}
 By Proposition \ref{GA1} and formula (\ref{equ-psi}), we see that  (\ref{eqinvd}) can be rewritten as
$$\mu(\Psi_1\cdots\Psi_{n-1})=\psi(\Psi_1)\cdots\psi(\Psi_{n-1}).  $$
Notice that for $i\geq 1$ and $i+2\leq k\leq n$, each term of $\Psi_i$ fixes $k$. Here we say an element $\sigma\in D_n$ fixes $k$ if and only if $\sigma$ maps $k$ to $k$. Thus $\Psi_i$ can be considered as an element of $\mathbb{Z}[D_j]$ for $i<j<n$.
Clearly, the weight function $\mu$ is also well-defined in this sense. Therefore we only need to show that
$$\mu(\Psi_1\cdots\Psi_{n-2}\,\Psi_{n-1})=\mu(\Psi_1\cdots\Psi_{n-2})\,\mu(\Psi_{n-1}).  $$
It suffices to prove that
\begin{equation}\label{equ-DD1}
  \mu(\sigma\cdot\Psi_{n-1})=\mu(\sigma)\cdot\mu(\Psi_{n-1})
\end{equation}
for any $\sigma=\sigma_1\cdots\sigma_{n-1}\in D_{n-1}$. Here $\sigma$ is considered as an element of $D_n$ which fixes $n$. It can be verified that
\begin{eqnarray*}
\sigma\cdot\Psi_{n-1}& = &
             \sigma_1\cdots\sigma_{n-1}n+\sigma_1\cdots\sigma_{n-2}n\sigma_{n-1}+\cdots+
\sigma_1n\cdots\sigma_{n-1}+n\sigma_1\cdots\sigma_{n-1}\\[5pt]
& &  +\bar{n}\bar{\sigma_1}\cdots\sigma_{n-1}+
\bar{\sigma_1}\bar{n}\cdots\sigma_{n-1}+\cdots+\bar{\sigma_1}\cdots\sigma_{n-1}\bar{n}.
\end{eqnarray*}
Then we have
 \begin{eqnarray*}
\lefteqn{\mu(\sigma\cdot\Psi_{n-1})} \\[5pt] & = &
             \mu(\sigma_1\cdots\sigma_{n-1}n)+\mu(\sigma_1\cdots\sigma_{n-2}n\sigma_{n-1})+\cdots+\psi(
\sigma_1n\cdots\sigma_{n-1})+\mu(n\sigma_1\cdots\sigma_{n-1}) \\[5pt]
  & & +\mu(\bar{n}\bar{\sigma_1}\cdots\sigma_{n-1})+\mu(
\bar{\sigma_1}\bar{n}\cdots\sigma_{n-1})+\cdots+
\mu(\bar{\sigma_1}\cdots\sigma_{n-1}\bar{n}) \\[5pt]
& = &\mu(\sigma)+qt\,\mu(\sigma)+\cdots+q^{n-2}t\,\mu(\sigma)+q^{n-1}t\,\mu(\sigma) \\[5pt]
 & &
 +q^{n-1}t\,\mu(\sigma)+q^n t\,\mu(\sigma)+\cdots+q^{2n-2} t\,\mu(\sigma)\\[5pt]
 & = & (1+t q^{n-1}+t q (1+q+\cdots+q^{2n-3}))\,\mu(\sigma).
\end{eqnarray*}
Hence (\ref{equ-DD1}) follows from (\ref{equ-psi}).
This completes the proof.  \qed

To prove formula (\ref{eqsord}) for  the bivariate generating function of $(sor_D,\tilde{l}_D')$, we recall   another factorization of the diagonal sum $\sum_{\sigma\in D_n}\sigma$ due to Petersen. For $2\leq j\leq n$, let
 \[ \Phi_j=1+\sum_{i\neq 0 \atop \bar{j}\leq i<j}t_{ij}.\]

\begin{prop}\label{GA2}
For $n\geq2$,
  $$\sum_{\sigma\in D_n}\sigma=\Phi_2\,\Phi_3\cdots\Phi_n. $$
\end{prop}

For an element $\sigma \in D_n$, we define  another weight function
$$ \nu(\sigma)=q^{sor_D(\sigma)}t^{\tilde{l}'_D(\sigma)}.$$
Again, the weight
function $\nu$ is considered as a linear map. It can be checked that
\begin{equation}\label{equ-phi}
  \nu(\Phi_i)=1+tq^{i-1}+tq(1+q+\cdots+q^{2i-3})=1+tq^{i-1}+tq[2i-2]_q.
\end{equation}

We conclude this paper with a proof of (\ref{eqsord}).

\noindent
{\em Proof of (\ref{eqsord}) in Theorem \ref{th4.2}.}
 By Proposition \ref{GA2} and (\ref{equ-phi}),  we find that
  (\ref{eqsord}) can be expressed in the following form
$$\nu(\Phi_2\cdots\Phi_n)=\nu(\Phi_2)\cdots\nu(\Phi_n).  $$
 As in the proof of (\ref{eqinvd}), we only need to show that
$$\nu(\Phi_2\cdots\Psi_{n})=\nu(\Phi_2\cdots\Phi_{n-1})\,\nu(\Phi_{n}).  $$
 It suffices to prove that
\begin{equation}\label{equ-DD2}
\nu(\sigma\cdot\Phi_n)=\nu(\sigma)\cdot\nu(\Phi_n),
\end{equation}
for any $\sigma=\sigma_1\cdots\sigma_{n-1}\in D_{n-1}$. Again, $\sigma$ is considered as an element of $D_n$ which fixes $n$.
 Since
\begin{eqnarray*}
\sigma\cdot\Phi_n& = &
             \sigma_1\cdots\sigma_{n-1}n + \sigma_1\cdots\sigma_{n-2}n\sigma_{n-1} + \cdots  +
\sigma_1n\cdots\sigma_{n-1}\sigma_2 + n\sigma_2\cdots\sigma_{n-1}\sigma_{1}\\[5pt]
& &  +\bar{n}\sigma_2\cdots\sigma_{n-1}\bar{\sigma_1}+
\sigma_1\bar{n}\cdots\sigma_{n-1}\bar{\sigma_2} + \cdots +
\bar{\sigma_1}\cdots\sigma_{n-1}\bar{n},
\end{eqnarray*}
we get
\begin{eqnarray*}
\lefteqn{\nu(\sigma\cdot\Phi_n)} \\[5pt] & = &
             \nu(\sigma_1\cdots\sigma_{n-1}n)+\nu(\sigma_1\cdots\sigma_{n-2}n\sigma_{n-1})+\cdots+\nu(
\sigma_1n\cdots\sigma_{n-1}\sigma_2)+\nu(n\sigma_2\cdots\sigma_{n-1}\sigma_1) \\[5pt]
  & & +\nu(\bar{n}\sigma_2\cdots\sigma_{n-1}\bar{\sigma_1})+
  \nu(\sigma_1\bar{n}\cdots\sigma_{n-1}\bar{\sigma_2})+\cdots+
\nu(\bar{\sigma_1}\sigma_2\cdots\sigma_{n-1}\bar{n}) \\[5pt]
& = &\nu(\sigma)+qt\,\nu(\sigma)+\cdots+q^{n-2}t\,\nu(\sigma)+q^{n-1}t\,\nu(\sigma) \\[5pt]
 & &
 +q^{n-1}t\,\nu(\sigma)+q^nt\,\nu(\sigma)+\cdots+q^{2n-2}t\,\nu(\sigma)\\[5pt]
 & = & (1+t q^{n-1}+t q (1+q+\cdots+q^{2n-3}))\,\nu(\sigma).
 \end{eqnarray*}
Hence (\ref{equ-DD2}) follows from (\ref{equ-phi}).
This completes the proof.   \qed

\vspace{0.5cm}
 \noindent{\bf Acknowledgments.}  This work was supported by  the 973
Project, the PCSIRT Project of the Ministry of Education,  and the National Science
Foundation of China.


\begin{thebibliography}{99}

\bibitem{BB} A. Bj\"{o}rner and F. Brenti, Combinatorics of Coxeter
Groups, Springer, New York, 2005.

\bibitem{Brenti} F. Brenti, $q$-Eulerian polynomials arising from Coxeter groups, European J. Combin. 15 (1994) 417-441.


\bibitem{Chen2} W.Y.C. Chen and R.P. Stanley, Derangements on the $n$-cube, Discrte Math. 115 (1993) 65-70.

\bibitem{Cori} R. Cori, Indecomposable permutations, hypermaps and
labeled Dyck paths, J. Combin. Theory Ser. A 116 (2009) 1326-1343.

\bibitem{OMR} P. Ossona de Mendez and P. Rosenstiehl, Transitivity and connectivity of permutations, Combinatorica 24 (2004) 487-502.

\bibitem{FH2} D. Foata and G.-H. Han, New permutation coding and
equidistribution of set-valued statistics, Theoret. Comput. Sci.
410 (2009) 3743-3750.

\bibitem{FH} D. Foata and G.-H. Han, Signed words and permutations, \Rmnum{1}:
a fundamental transformation, Proc. Amer. Math. Soc. 135 (2007) 31-40.


\bibitem{knuth} D.E. Knuth, The Art of Computer Programming,
Vol. 3, Addison-Wesley, 1998.


\bibitem{Lehmer} D.H. Lehmer, Teaching combinatorial tricks to a computer, in: Proc. Sympos. Appl. Math., Vol. 10, Amer. Math. Soc, Providence, RI, 1960, pp. 179-193.

\bibitem{Petersen} T.K. Petersen, The sorting index, Adv. Appl. Math. 47 (2011) 615-630.

\bibitem{Wilson1} M.C. Wilson, An interesting new Mahonian permutation statistic, Electron. J. Combin. 17 (2010) R147.

\bibitem{Wilson2} M.C. Wilson, Random and exhaustive generation of permutations and cycles, Ann. Combin. 12 (2009) 509-520.

\end{thebibliography}
\end{document}